\documentclass[reqno,12pt]{amsart}
\usepackage{a4wide}
\usepackage [frenchb]{babel}
\numberwithin{equation}{section}

\newtheorem{Theorem}{Th{\'e}or{\`e}me}

\newtheorem{Lemma}{Lemme}
\theoremstyle{remark}

\def\al{\alpha}

\def\({\left(}
\def\){\right)}
\def\[{\left[}
\def\]{\right]}
\def\ii{\infty}

\def\li{\operatorname{Li}}
\def\dis{\displaystyle}
\def\dd{\textup{d}}

\begin{document}

\title[]{Approximants de Pad{\'e} des $q$-polylogarithmes}
\author[]{C. Krattenthaler$^\dagger$
et T. Rivoal}

\address{
Institut Girard Desargues,
Universit{\'e} Claude Ber\-nard Lyon-I,
21, avenue Claude Ber\-nard,
F-69622 Villeurbanne Cedex, France}
\email{kratt@euler.univ-lyon1.fr}
\address{
Laboratoire de Math{\'e}matiques Nicolas Oresme, CNRS UMR 6139,
Universit{\'e} de Caen,  BP 5186,
14032 Caen cedex,
France}
\email{rivoal@math.unicaen.fr}
\thanks{$^\dagger$ Recherche partiellement support{\'e}e
par le Programme {\og Accro{\^\i}tre le potentiel humain de recherche \fg}
de la Commission Europ{\'e}enne, contrat HPRN-CT-2001-00272,
``Algebraic Combinatorics in Europe"}

\dedicatory{D{\'e}di{\'e} {\`a} Wolfgang Schmidt, pour son soixante-dixi\`eme anniversaire}



\begin{abstract}
Nous r{\'e}solvons un probl{\`e}me d'approximation simultan{\'e}e de type  
Pad{\'e} mettant en jeu des $q$-analogues sp{\'e}cifiques
des polylogarithmes et des puissances du logarithme. Ce probl{\`e}me
est fortement li{\'e} aux
r{\'e}sultats r{\'e}cents des auteurs et de Wadim Zudilin
[{\og S{\'e}ries hyperg{\'e}om{\'e}triques basiques, fonction
$q$-z{\^e}ta et s{\'e}ries d'Eisenstein\fg},
{\`a} para{\^\i}tre au J. Inst.\ Math.\ Jussieu] sur la dimension
de l'espace vectoriel engendr{\'e} par des $q$-analogues des
valeurs de la fonction z{\^e}ta de Riemann aux nombres entiers.
Nous montrons aussi que ce  r{\'e}sultat peut {\^e}tre consid\'er{\'e}
comme $q$-analogue d'un r{\'e}sultat de St\'ephane Fischler et du deuxi{\`e}me auteur
[J. Math.\ Pures  Appl.\ {\bf 82} (2003), 1369--1394].

\medskip

\noindent
{\sc Abstract.}
We solve a Pad{\'e}-type problem of approximating three specific
functions simultaneously by $q$-analogues
of polylogarithms,
respectively by powers of the logarithm. This problem is
intimately related to recent results of the authors and Wadim Zudilin
[{\og S{\'e}ries hyperg{\'e}om{\'e}triques basiques, fonction
$q$-z{\^e}ta et s{\'e}ries d'Eisenstein\fg},
J. Inst.\ Math.\ Jussieu (to appear)] on the dimension
of the vector space generated by $q$-analogues of
values of the Riemann zeta function at integers.
We also show that our result can be considered as a $q$-analogue
of a result of~St\'ephane Fischler and the second author
[J. Math.\ Pures  Appl.\ {\bf 82} (2003), 1369--1394].

\end{abstract}

\maketitle

%%%%%%%%%%%%%%%%%%%%%%%%%%%%%%%%%%%
%%%%%%%%%%%%%%%%%%%%%%%%%%%%%%%%%%
%%%%%%%%%%%%%%%Introduction%%%%%%
%%%%%%%%%%%%%%%%%%%%%%%%%%%%%%%%%%
%%%%%%%%%%%%%%%%%%%%%%%%%%%%%%%%%

\section{Introduction}

Consid{\'e}rons la s{\'e}rie
$$
\zeta_q(s)=\sum_{k= 1}^{\infty}k^{s-1}\frac{q^{k}}{1-q^k},
$$
qui converge pour tout complexe $\vert q \vert <1$ et tout entier $s\ge 1$.
La notation $\zeta_q$ est justifi{\'e}e par le fait
que cette fonction est un $q$-analogue
de la fonction z{\^e}ta de Riemann $\zeta(s)$ au sens suivant
(voir \cite[paragraphe~4.1]{krz}, \cite[Theorem~2]{jap} ou \cite{zud5}),
\begin{equation*}
\lim_{q\to 1}\, (1-q)^s\zeta_q(s)=(s-1)!\sum_{k=1}^{\infty}\frac{1}{k^s}
=(s-1)!\,\zeta(s).
\end{equation*}
Dans \cite{krz}, les deux auteurs et W. Zudilin ont montr{\'e} que la dimension de l'espace
vectoriel engendr{\'e} sur $\mathbb{Q}$ par 1, $\zeta_q(3)$, $\zeta_q(5)$, ...,
$\zeta_q(A)$ ($A\ge 3$ impair) est minor{\'e}e par
$\dis \frac{\pi+o(1)}{2\sqrt{\pi^2+12}} \,\sqrt{A}$ lorsque
$1/q\in\mathbb{Z}\setminus\{\pm 1\}$.
La d{\'e}monstration utilise les fonctions $q$-polylogarithmes, d{\'e}finies
pour tout entier $s\ge 1$, par
\begin{equation}\label{eq:Li}
\li_s(z;q)=\sum_{k= 1}^{\infty}\frac{q^k}{(1-q^k)^s}\,z^{k},
\end{equation}
o{\`u} $z$ et $q$ d{\'e}signent des nombres
complexes tels que
$\vert q\vert <1$ et $\vert zq\vert <1$. Ces fonctions constituent des $q$-analogues des
polylogarithmes usuels $\li_j(z)$ au sens  suivant~:
$$
\lim_{q\to 1}\, (1-q)^s\li_s(z;q)=\sum_{k=1}^{\infty}\frac{z^{k}}{k^s}=\li_s(z).
$$
Notons que les polylogarithmes  sont aussi utilis{\'e}s au cours de la d{\'e}montration du
th{\'e}or{\`e}me suivant (dont celui rappel{\'e} ci-dessus est
un $q$-analogue)~: la dimension de l'espace vectoriel engendr{\'e} sur $\mathbb{Q}$
par 1, $\zeta(3)$, $\zeta(5)$, ..., $\zeta(A)$ ($A\ge 3$ impair) est
minor{\'e}e par $\dis \frac{1+o(1)}{1+\log(2)} \, \log(A)$ (voir ~\cite{br, ri}).
Dans les deux cas, la d{\'e}monstration est en fait bas{\'e}e sur
une {\'e}tude tr{\`e}s fine d'une s{\'e}rie
($q$-)hyperg{\'e}om{\'e}trique bien choisie
que l'on commence
par exprimer comme une combinaison lin{\'e}aire polynomiale en les
($q$-)polylogarithmes. Plus pr{\'e}cis{\'e}ment, soient $A, n, r$ des entiers positifs tels que
$0\le r\le A/2$. D{\'e}finissons les factorielles d{\'e}cal{\'e}es
(ou symboles de Pochhammer)
$(\al)_m=\al(\al+1)\cdots(\al+m-1)$ et
les factorielles $q$-d{\'e}cal{\'e}es
$(\al;q)_m=(1-\al)(1-\al q)\cdots(1-\al q^{m-1})$, avec la convention
usuelle que les produits vides pour $m=0$ valent 1. On pose alors
$$
S_n(z;q)= (q;q)_{n}^{A-2r}
\sum_{k=1}^{\ii}q^k\frac{(q^{k-rn};q)_{rn}\,(q^{k+n+1};q)_{rn}}
{(q^{k};q)_{n+1}^A}\,q^{(k-1/2)(A-2r)n/2}z^{-k},
$$
avec $\vert q\vert <1 \le \vert z\vert $ et $A$ pair, ainsi que
$$
S_n(z)=n!^{A-2r}\sum_{k=1}^{\infty}
\frac{(k-rn)_{rn}\,(k+n+1)_{rn}}{(k)_{n+1}^{A}}\,z^{-k},
$$
avec $\vert z\vert \ge 1$.
Il est alors facile de montrer l'existence de deux familles de polyn{\^o}mes
$P_{j,n}(z;q)\in\mathbb{C}(q)[z]$ et $P_{j,n}(z)\in\mathbb{C}[z]$,
de degr{\'e} au plus $n$, tels que
\begin{equation}\label{eq:qlineaireforme}
S_{n}(z;q)=P_{0,n}(z;q)+\sum_{j=1}^A P_{j,n}(z;q)\li_j(1/z;q)
\end{equation}
et
\begin{equation}\label{eq:lineaireforme}
S_{n}(z)=P_{0,n}(z)+\sum_{j=1}^A P_{j,n}(z)\li_j(1/z).
\end{equation}
Par ailleurs, les
num{\'e}rateurs des sommandes dans les d{\'e}finitions de $S_n(z;q)$ et $S_n(z)$ s'annulent pour
les indices $k\in\{1, \ldots, rn\}$, ce qui assure que
l'ordre en $z=0$ des deux s{\'e}ries est exactement $rn+1$~: les
{\'e}quations~\eqref{eq:qlineaireforme} et~\eqref{eq:lineaireforme}
peuvent donc {\^e}tre vues comme des probl{\`e}mes d'approximations de type
Pad{\'e} pour les s{\'e}ries enti{\`e}res 1 et $\li_j(z;q)$, respectivement 1 et
$\li_j(z)$. L'information n'est cependant pas suffisante pour affirmer
qu'il n'existe, {\`a} constante multiplicative pr{\`e}s, qu'une seule fonction $S_n(z;q)$,
resp. $S_n(z)$, v{\'e}rifiant~\eqref{eq:qlineaireforme},
resp.~\eqref{eq:lineaireforme}, et qui s'annule {\`a} l'ordre $rn+1$.

Dans~\cite{firi}, sont {\'e}nonc{\'e}es des conditions suppl{\'e}mentaires, de type Pad{\'e}, portant
sur des objets  li{\'e}s {\`a}
la s{\'e}rie $S_n(z)$ et qui suffisent {\`a} assurer que $S_n(z)$ est bien la seule solution
de~\eqref{eq:lineaireforme} pour des polyn{\^o}mes $P_{j,n}(z)$ de degr{\'e} au plus $n$. Voici
l'{\'e}nonc{\'e} pr{\'e}cis.

{\it \'Etant donn\'es des entiers $A\ge 1$,  $n\ge 0$, $\rho\ge 0$,
$\sigma\ge 0$  tels que $\rho+\sigma+2\le A(n+1)$,
on cherche {\`a} r{\'e}soudre le probl{\`e}me d'approximations simultan{\'e}es
de Pad{\'e} suivant~: \it d{\'e}terminer des polyn{\^o}mes
{\em(}d{\'e}pendants  de
$A$, $n$, $\rho$, $\sigma${\em)}  $P_{0,n}(z)$,
$\overline{P}_{0}(z)$ et
$P_{j}(z)$ {\em(}pour $j=1,\ldots,A${\em)}, de degr{\'e} au plus $n$
et {\`a} coefficients dans $\mathbb{Q}$, tels que }
\begin{equation}\label{problemepade}
\begin{cases}
\dis S(z)=P_{0}(z)+\sum_{j=1}^A P_{j}(z)\li_j(1/z)=
\mathcal{O}(z^{-\rho-1}) \quad\textup{quand}\quad  z\to \infty \, ;\\
\dis \overline{S}(z)=\overline{P}_{0}(z;q)+
\sum_{j=1}^A P_{j}(z)\li_j(z)=\mathcal{O}(z^{\sigma+n+1})
\quad \textup{quand}\quad z\to 0 \,;\\
\dis I(z)=\sum_{j=1}^A P_{j}(z)\frac{\log^{j-1}(1/z)}{(j-1)!}
=\mathcal{O}((z-1)^{A(n+1)-\rho-\sigma-2})\quad \textup{quand}\quad z\to 1.
\end{cases}
\end{equation}
(Ici et dans toute la suite, la fonction logarithme est d{\'e}finie avec sa
branche principale~: $\log(z)=\log\vert z\vert+i\,\textup{arg}(z)$, avec
$-\pi<\textup{arg}(z)\le \pi$. On notera $\mathbb{R}_{-}$ l'ensemble
des r{\'e}els n{\'e}gatifs.) On a alors le r{\'e}sultat suivant, qui
r{\'e}sout conpl{\`e}tement ce probl{\`e}me.
\begin{Theorem}\label{thm:firi}Dans les conditions ci-dessus,
le probl{\`e}me~\eqref{problemepade} a une solution unique,
 {\`a} une constante multiplicative pr{\`e}s. En choisissant
cette constante {\'e}gale {\`a} 1, on a
$$
S(z)=\sum_{k=1}^{\ii}\frac{(k-\rho)_{\rho}\,(k+n+1)_{\sigma}}
{(k)_{n+1}^A}\,z^{-k} \quad \textup{et} \quad
I(z)=\frac{1}{2i\pi}\int_{\mathcal{C}}
\frac{(s-\rho)_{\rho}\,(s+n+1)_{\sigma}}
{(s)_{n+1}^A} \,z^{-s} \,\dd s,
$$
o{\`u} $\mathcal{C}$ est n'importe quelle
courbe ferm{\'e}e orient{\'e}e dans le sens direct
qui entoure les p{\^o}les de
l'int{\'e}grande, {\em i.e.}
$0, -1, \ldots, -n$.
\end{Theorem}

\medskip

Le but de cette note est de prouver un $q$-analogue  du
th{\'e}or{\`e}me pr{\'e}c{\'e}dent. On suppose dor{\'e}navant que $q\not\in
\mathbb{R}_{-}$ ce qui permet de d{\'e}finir $\log_q(z)=\log(z)/\log(q)$.

{\it \'Etant donn\'es des entiers $A\ge 1$,  $n\ge 0$, $\rho\ge 0$,
$\sigma\ge 0$ et $\nu\ge 0$ tels que $\rho+\sigma+\nu+2\le A(n+1)$,
on cherche {\`a} r{\'e}soudre le probl{\`e}me d'approximations simultan{\'e}es
de Pad{\'e} suivant~: d{\'e}terminer des polyn{\^o}mes
{\em(}d{\'e}pendants  de
$A$, $n$, $\rho$, $\sigma$ et $\nu${\em)}  $P_{0}(z;q)$,
$\overline{P}_{0}(z;q)$ et
$P_{j}(z;q)$ {\em(}pour $j=1,\ldots,A${\em)} en la variable $z$, de degr{\'e} au plus $n$
et {\`a} coefficients dans $\mathbb{Q}(q)$, tels que }

\begin{equation}\label{problemepade1}
\begin{cases}
\dis S(z;q)=P_{0}(z;q)+\sum_{j=1}^A P_{j}(z;q)\li_j(1/z;q)=
\mathcal{O}(z^{-\rho-1}) \quad\textup{quand}\quad  z\to \infty \, ;\\
\dis \overline{S}(z;q)=\overline{P}_{0}(z;q)+
\sum_{j=1}^A P_{j}(z;q)\li_j(z;1/q)=\mathcal{O}(z^{\sigma+n+1})
\quad \textup{quand}\quad z\to 0 \,;\\
\dis I(z;q)=-\sum_{j=1}^A
P_{j}(zq^{1-j};q)\frac{\(-\log_q(1/z)\)_{j-1}}{(j-1)!}
=\mathcal{O}(z-q^{-\ell})\quad \textup{quand}\quad z\to q^{-\ell}\\
\qquad\quad
\textup{pour tout} \;\ell\in\{-\nu, -\nu+1,  \ldots, A(n+1)-\rho-\sigma-\nu-2\}.
\end{cases}
\end{equation}
\begin{Theorem} \label{thm:1}
Dans les conditions ci-dessus, le probl{\`e}me~\eqref{problemepade1} a une
solution unique, {\`a} une constante multiplicative pr{\`e}s. En choisissant
cette constante {\'e}gale {\`a} 1, on a
$$
S(z;q)=\sum_{k=1}^{\ii}q^{k}\frac{(q^{k-\rho};q)_{\rho}\,(q^{k+n+1};q)_{\sigma}}
{(q^{k};q)_{n+1}^A}\,q^{\nu k}z^{-k}
$$
et
$$
I(z;q)=\frac{1}{2i\pi}\int_{\mathcal{C}}
\frac{(sq^{-\rho};q)_{\rho}\,(sq^{n+1};q)_{\sigma}}
{(s;q)_{n+1}^A}s^{\nu-\log_q(z)} \,\dd s,
$$
o{\`u} $\mathcal{C}$ est n'importe quelle courbe ferm{\'e}e orient{\'e}e dans le sens direct
qui entoure les p{\^o}les de
l'int{\'e}grande, {\em i.e.}
$1, q^{-1}, \ldots, q^{-n}$, sans traverser la coupure $\mathbb{R}_{-}$.
\end{Theorem}

Avant de passer {\`a} la d{\'e}monstration du Th{\'e}or{\`e}me~\ref{thm:1},
faisons quelques remarques~:

Les  Th{\'e}or{\`e}mes~\ref{thm:firi} et~\ref{thm:1} sont formellement
tr{\`e}s similaires, {\`a} ceci pr{\`e}s que le param{\`e}tre {\og
  $q$-analogique\fg} $\nu$ n'a pas d'{\'e}quivalent dans le cas
classique. La diff{\'e}rence majeure se situe dans l'{\'e}nonc{\'e} des
conditions d'annulation des fonctions $I(z)$ et $I(z;q)$~: cela
n'a cependant rien surprenant, puisqu'il est fr{\'e}quent dans ce genre de
situation que des singularit{\'e}s en certaines puissances de $q$
confluent vers une unique singularit{\'e} en 1
(avec une certaine multiplicit\'e)
lorsque $q\to 1$. Nous explicitons plus en d{\'e}tail la {\og convergence\fg} du
Th{\'e}or{\`e}me~\ref{thm:1} vers le Th{\'e}or{\`e}me~\ref{thm:firi}
 au paragraphe~\ref{sec:confl}.

Il est {\`a} noter l'utilisation, naturelle
dans notre contexte, de $\log(z)/\log(q)$ com\-me
$q$-ana\-lo\-gue de la fonction logarithme. Cet analogue, qui poss{\`e}de donc
une monodromie non-triviale en 0, est un choix historiquement classique~:
voir~\cite{adams}.
Certaines th{\'e}ories g{\'e}om{\'e}triques
r{\'e}centes ({\'e}tudiant l'analogie entre {\'e}quations aux
$q$-diff{\'e}rences et {\'e}quations diff{\'e}rentielles) ont mis en avant
un $q$-analogue diff{\'e}rent du logarithme~:
J. Sauloy~\cite{sau} utilise comme $q$-logarithme la
fonction $\ell_q(z)=z\theta^{\prime}_q(z)/\theta_q(z)$, avec
$\theta_q(z)=\sum_{n\in\mathbb{Z}} (-1)^n q^{-n(n-1)/2}z^n$, qui est
m{\'e}romorphe sur $\mathbb{C}$ et dont les  p{\^o}les
confluent le long d'une spirale lorsque $q\to 1$, cette spirale
agissant alors comme une coupure du plan pour le logarithme usuel.

\section{D{\'e}monstration du Th{\'e}or{\`e}me \ref{thm:1}}

Remarquons tout d'abord que les polyn{\^o}mes $P_{0}(z;q)$ et $\overline
P_{0}(z;q)$ sont d{\'e}termin{\'e}s de fa{\c c}on unique, une fois connus les autres
polyn{\^o}mes du probl{\`e}me~\eqref{problemepade1}. De plus, celui-ci se traduit par un
syst{\`e}me lin{\'e}aire dont les inconnues sont les $A(n+1)$ coefficients
des polyn{\^o}mes $P_{j}(z;q)$ ($j\ge 1$) et dont le nombre
d'{\'e}quations est  $A(n+1)-1$~: il y a donc au moins une solution
non-triviale ({\it i.e.} non identiquement nulle).

On utilise temporairement la notation $P_{j}(z;q)$ pour d{\'e}signer des
polyn{\^o}mes g{\'e}n{\'e}riques
de degr{\'e} au plus $n$ en $z$ et {\`a} coefficients dans $\mathbb{Q}(q)$,
sans pr{\'e}sager qu'il s'agisse
des solutions du probl{\`e}me \eqref{problemepade1}. On pose
\begin{equation} \label{eq:P_j}
P_{j}(z;q)=\sum_{t=0}^n p_{j,t}(q)z^t,
\end{equation}
o{\`u} $p_{j,t}(q)\in \mathbb{Q}(q)$, de telle sorte que
$$
\sum_{j=1}^A P_{j}(z;q)\li_j(1/z;q)=\sum_{k=1-n}^{\ii}q^kz^{-k}
\sum_{j=1}^A\sum_{t=\max(0,1-k)}^n\frac{q^tp_{j,t}(q)}{(1-q^{k+t})^j}.
$$
Il est utile {\`a} ce point d'introduire la fraction rationnelle (qui
d{\'e}pend aussi de A)~:
\begin{align} \label{eq:R1}
R(s;q)&=\sum_{j=1}^A\sum_{t=0}^n\frac{q^t p_{j,t}(q)}{(1-sq^{t})^j}\\
\label{eq:R2}
&=\frac{\Pi(s;q)}{(s;q)_{n+1}^A},
\end{align}
o{\`u} $\Pi(s;q)$ est un polyn{\^o}me en $s$, {\`a} coefficients dans
$\mathbb{Q}(q)$,
et de degr{\'e} $<A(n+1)$ puisque la d{\'e}composition en
{\'e}l{\'e}ments simples~\eqref{eq:R1} de $R(s;q)$ est sans
partie principale. Il est clair que la connaissance de  $\Pi(s;q)$
d{\'e}termine  {\it de facto} les polyn{\^o}mes $P_{j}(z;q)$
($j\ge 1$). On en d{\'e}duit que si l'on connait $\Pi(s;q)$, alors on a
$$
S(z;q)=\sum_{k=1}^{\ii} q^kR(q^k;q) z^{-k}
$$
et aussi, par un simple calcul de r{\'e}sidus
utilisant l'expression~\eqref{eq:R1},
$$
I(z;q)=\frac{1}{2i\pi}\int_{\mathcal{C}}
R(s;q)s^{-\log_q(z)} \,\dd s,
$$
o{\`u} $\mathcal{C}$ est n'importe quelle
courbe ferm{\'e}e entourant dans le sens direct les p{\^o}les de $R(z;q)$
et qui ne traverse pas la coupure $\mathbb{R}_{-}$. Nous allons
maintenant montrer que la solution du
probl{\`e}me de Pad{\'e}~\eqref{problemepade1} est unique, {\`a} une constante
multiplicative pr{\`e}s, et la d{\'e}terminer en
explicitant le polyn{\^o}me {\og codant \fg} $\Pi(s;q)$. Pour cela,
nous interpr{\'e}tons chacune des  conditions
de~\eqref{problemepade1} par quatre lemmes~: il en d{\'e}coulera alors
que
$$
\Pi(s;q)=(sq^{-\rho};q)_{\rho}\,(sq^{n+1};q)_{\sigma}\,s^{\nu},
$$
{\`a} une constante multiplicative pr{\`e}s.

\begin{Lemma} \label{lem:1}
Les polyn{\^o}mes $P_{1}, \ldots, P_{A}$ v{\'e}rifient la
  premi{\`e}re condition de~\eqref{problemepade1} si, et seulement si,
$$
\prod_{i=1}^{\rho} (1-sq^{-i})=(sq^{-\rho};q)_{\rho}
\quad\textup{divise}\quad  \Pi(s;q).
$$
\end{Lemma}

\begin{proof}
La premi{\`e}re condition de~\eqref{problemepade1} se traduit par
l'annulation des coefficients de Taylor de $S(z;q)$ d'indices $1, 2, \ldots, \rho$, ce qui
{\'e}quivaut {\`a}
l'annulation de la fonction
$k\mapsto \Pi(q^k;q)$ en $k=1, 2, \ldots, \rho$. Cela {\'e}quivaut en fait {\`a}
l'annulation du polyn{\^o}me $\Pi(s;q)$ en
$s= q, q^2, \ldots,
q^{\rho}$, d'o{\`u} l'assertion.
\end{proof}

\begin{Lemma}\label{lem:2}
Les polyn{\^o}mes $P_{1}, \ldots, P_{A}$ v{\'e}rifient la
deuxi{\`e}me condition de~\eqref{problemepade1} si, et seulement si,
$$
\prod_{i=n+1}^{n+\sigma} (1-sq^{i})=(sq^{n+1};q)_{\sigma}
\quad \textup{divise}\quad  \Pi(s;q).
$$
\end{Lemma}

\begin{proof} Dans la deuxi{\`e}me condition de~\eqref{problemepade1},
  on change $z$ en $1/z$, puis on multiplie par $z^n$, de telle sorte
  que $z^n\overline{S}(1/z;q)=\mathcal{O}(z^{-\sigma-1})$. On a
\begin{eqnarray*}
z^n\sum_{j=1}^A P_{j}(1/z;q)\li_j(1/z;1/q) &=&
\sum_{k=1-n}^{\ii}q^{-k}z^{-k}
\sum_{j=1}^A\sum_{t=\max(0,1-k)}^n\frac{q^{-t}p_{j,n-t}(q)}{(1-q^{-k-t})^j}\\
&=& q^{-n} \sum_{k=1-n}^{\ii}q^{-k}z^{-k}
\sum_{j=1}^A\sum_{t=0}^{\min(n,n+k-1)}
\frac{q^tp_{j,t}(q)}{(1-q^{t-k-n})^j}.
\end{eqnarray*}
Pour $k\ge 1$, le coefficient de $z^{-k}$ dans la s{\'e}rie
$z^n\overline S(1/z;q)$ est donc donn{\'e} par\break
$q^{-n-k}R(q^{-k-n};q)$. L'annulation des coefficients
de Taylor de
$\overline S(z;q)$ d'indices $1, 2, \ldots, \sigma$ {\'e}quivaut
donc {\`a} l'annulation du polyn{\^o}me $\Pi(s;q)$ en $s=q^{-n-1},  \ldots,
q^{-n-\sigma}$, d'o{\`u} l'assertion.
\end{proof}

\begin{Lemma}
\label{lem:3} La condition $I(q^{-j};q)=0$ pour $j\in\{-\nu, \ldots,-1\}$
{\'e}quivaut {\`a}
$$
s^{\nu} \quad \textup{divise}\quad  \Pi(s;q).
$$
\end{Lemma}

\begin{proof}
Rappelons que
$$
I(z;q)=\frac{1}{2i\pi}\int_{\mathcal{C}}
 \frac{\Pi(s;q)}{(s;q)_{n+1}^{A}} s^{-\log_q(z)} \,\dd s.
$$
o{\`u} $\mathcal{C}$ est n'importe quelle courbe ferm{\'e}e orient{\'e}e dans le
sens direct
et entourant les p{\^o}les de
l'int{\'e}grande, {\it i.e.}
1, $q^{-1}, \ldots, q^{-n}$, sans traverser la coupure
$\mathbb{R}_{-}$.

Soit $-n\le j\le -1$, et supposons que $I(q^{-j},q)=0$.
Alors, pour tout contour $\mathcal{C}$ entourant
les points $1, q^{-1}, \ldots, q^{-n}$ mais ne traversant
pas la coupure, donc n'entourant pas 0, on a
$$
0=I(q^{-j};q)=\frac{1}{2i\pi}\int_{\mathcal{C}} \frac{
\Pi(s;q)}{s^{-j}(s;q)_{n+1}^{A}}  \,\dd s,
$$
o{\`u} l'int{\'e}grande est une fraction rationnelle dont 0 est peut-{\^e}tre un
p{\^o}le~: nous allons montrer que ce n'est ({\'e}ventuellement) le cas
que si $j<-\nu$, ce qui prouvera que $s^\nu$
divise   $\Pi(s;q)$.

Pour cela, notons
que, par le th{\'e}or{\`e}me des r{\'e}sidus  et parce que le r{\'e}sidu {\`a}
l'infini de l'int{\'e}grande $F(s;q)$ de $I(q^{-j};q)$  est nul (car le degr{\'e} du
num{\'e}rateur de $F(s;q)$ est $\le A(n+1)-1$,
et celui de son d{\'e}nominateur
$\ge A(n+1)+1$), on a
\begin{multline*}
0=-\textup{Res}_{\ii}(F)=\frac{1}{2i\pi}\int_{\mathcal{C}'} F(s;q)\,\dd
s \\ = \frac{1}{2i\pi}\int_{\mathcal{C}} F(s;q)  \,\dd s +
\frac{1}{2i\pi}\int_{\mathcal{C}_0} F(s;q)  \,\dd s=
\frac{1}{2i\pi}\int_{\mathcal{C}_0} F(s;q)  \,\dd s,
\end{multline*}
o{\`u} $\mathcal{C}'$ est un cercle entourant tous les p{\^o}les de $F$ ainsi que 0, et
$\mathcal{C}_0$ un cercle  entourant 0 et aucun autre p{\^o}le de
$F$, les deux orient{\'e}s dans le sens direct. Donc
\begin{equation}\label{eq:annulationOmega}
0=\frac{1}{2i\pi}\int_{\mathcal{C}_0} F(s;q)  \,\dd s=\frac{1}{(-j-1)!}
\((s;q)_{n+1}^{-A}\Pi(s;q)\)^{(-j-1)}\bigg\vert_{s=0}.
\end{equation}
Puisque la fonction $s\mapsto (s;q)_{n+1}^{-A}$ ne s'annule pas en $s=0$,
on d{\'e}duit de~\eqref{eq:annulationOmega}, par
r{\'e}currence sur $j\in\{-1,-2, \ldots,- \nu\}$, que
$0=\Pi(0;q)=\Pi^{(1)}(0;q)=\cdots=\Pi^{(\nu-1)}(0;q)$, ce qui prouve que $s^{\nu}$
divise $\Pi(s;q)$. La r{\'e}ciproque se montre facilement en renversant cet argument.
\end{proof}

\begin{Lemma}
La condition $I(q^{-j};q)=0$ pour $j\in\{0, \ldots, A(n+1)-\rho-\sigma-\nu-2\}$
{\'e}quivaut {\`a}
$$
\textup{deg}(\Pi)\le \rho+\sigma+\nu.
$$
\end{Lemma}
\begin{proof}
D{\'e}veloppons $\Pi(s;q)/(s;q)_{n+1}^A$ en s{\'e}rie enti{\`e}re en
$s=\infty$~:
$$
\frac{\Pi(s;q)}{(s;q)_{n+1}^A}=\sum_{k=\omega}^{\infty} \frac{c_k}{s^k}
$$
o{\`u} $\omega=A(n+1)-\textup{deg}(\Pi)$ est l'ordre de cette fraction
rationnelle {\`a} l'infini, et les $c_k$ sont des nombres
complexes. Notons que cette s{\'e}rie converge au moins pour $\vert s\vert$ assez
grand puisque $\omega \ge 1$, disons pour $\vert s\vert\ge S$.
On choisit alors un cercle suffisamment grand
pour que l'on puisse int{\'e}grer cette s{\'e}rie terme {\`a} terme, disons
le cercle $\overline{\mathcal{C}}=\{z:\vert z\vert=S+1\}$.
On a alors
$$
I(q^{-j};q)=\frac{1}{2i\pi} \int_{\overline{\mathcal{C}}}
\frac{\Pi(s;q)}{(s;q)_{n+1}^A} s^j \,\dd s = \sum_{k=\omega}^{\infty} c_k
\int_{\overline{\mathcal{C}}} s^{j-k}\,\dd s =c_{j+1}.
$$
L'annulation de $I(q^{-j};q)$ pour $j\in\{0, \ldots,
A(n+1)-\rho-\sigma-\nu-2\}$ {\'e}quivaut donc
{\`a}
$$
\omega \ge A(n+1)-\rho-\sigma-\nu,
$$
ce qui {\'e}quivaut plus simplement {\`a} $\textup{deg}(\Pi)\le \rho+\sigma+\nu$.
\end{proof}

\begin{proof}[D{\'e}monstration du Th{\'e}or{\`e}me~\ref{thm:1}] Puisque les
  polyn{\^o}mes $s^\nu$, $(sq^{-\rho};q)_{\rho}$ et
  $(sq^{n+1};q)_{\sigma}$ n'ont pas de racines communes, les trois
  premiers lemmes montrent que
  $s^\nu(sq^{-\rho};q)_{\rho}\,(sq^{n+1};q)_{\sigma}$ divise
  $\Pi(s;q)$. Or le dernier lemme montre que $\textup{deg}(\Pi)\le
  \rho+\sigma+\nu$, ce qui ach{\`e}ve la d{\'e}monstration.
\end{proof}

\section{Confluence du Th{\'e}or{\`e}me~\ref{thm:1} vers le  Th{\'e}or{\`e}me~\ref{thm:firi}}
\label{sec:confl}

Dans ce paragraphe, nous explicitons le sens pr{\'e}cis en lequel
le Th{\'e}or{\`e}me~\ref{thm:1} {\og tend\fg} vers le
Th{\'e}or{\`e}me~\ref{thm:firi}.

Pour commencer, on remarque que, {\'e}videmment, on a
$$
\lim _{q\to1} ^{}(1-q)^{A(n+1)-\sigma-\rho}S(z;q)=S(z)\quad
\text{et}\quad
\lim _{q\to1} ^{}(1-q)^{A(n+1)-\sigma-\rho}\overline S(z;q)=\overline S(z).$$
De plus, en faisant la substitution $s\to q^t$ dans l'int{\'e}grale
d{\'e}finissant $I(z;q)$, on voit que
$$
\lim _{q\to1} ^{}(1-q)^{A(n+1)-\sigma-\rho-1}I(z;q)=-I(q).
$$
D'autre part, la d{\'e}finition~\eqref{eq:P_j} des polyn{\^o}mes
$P_j(z;q)$ est donn{\'e}e par leurs coefficients $p_{j,t}(q)$ qui
figurent dans \eqref{eq:R1}, avec
$\Pi(s;q)=(sq^{-\rho};q)_{\rho}\,(sq^{n+1};q)_{\sigma}\,s^{\nu}$.
Explicitement, les $p_{j,t}(q)$ sont donn{\'e}s par
$$
p_{j,t}(q)=\frac {(-1)^{A-j}q^{t(A-j-1)}} {(A-j)!}
\frac {\partial^{A-j}} {\partial s^{A-j}}\left(\frac {(1-sq^t)^{A}\Pi(s;q)}
{(s;q)^A_{n+1}}\right)\bigg\vert_{s=1}.
$$
Par cons{\'e}quent, $p_{j,t}(q)$ est une fraction rationnelle en $q$~;
si cette fraction rationnelle est {\'e}crite sous forme
r{\'e}duite, la plus grande puissance de $1-q$ qui divise le
d{\'e}nominateur est
$(1-q)^{A(n+1)-\sigma-\rho-j}$.
En particulier, la limite $\lim _{q\to1}
^{}(1-q)^{A(n+1)-\sigma-\rho-j}P_{j}(z;q)$
existe~: c'est un  polyn{\^o}me en $z$, que l'on notera $Q_j(z)$.
De fa\c{c}on similaire, 
la limite $\lim _{q\to1} ^{}(1-q)^{A(n+1)-\sigma-\rho}\overline
P_{0}(z;q)$ existe et sera not{\'e}e $\overline Q_0(z)$.

Si l'on combine ces remarques avec \eqref{eq:Li} et le fait que
$$
\lim _{q\to1} ^{}(1-q)\log_q(z)=-\log(z),
$$
on en d{\'e}duit que, en multipliant les deux premi{\`e}res conditions
dans \eqref{problemepade1} par\break
$(1-q)^{A(n+1)-\sigma-\rho}$, en multipliant
la troisi{\`e}me par $(1-q)^{A(n+1)-\sigma-\rho-1}$, puis en faisant
finalement tendre $q$ vers 1, on obtient
\begin{equation}\label{problemepadeQ}
\begin{cases}
\dis S(z)=Q_{0}(z)+\sum_{j=1}^A Q_{j}(z)\li_j(1/z)=
\mathcal{O}(z^{-\rho-1}) \quad\textup{quand}\quad  z\to \infty \, ;\\
\dis \overline{S}(z)=\overline{Q}_{0}(z)+
\sum_{j=1}^A Q_{j}(z)\li_j(z)=\mathcal{O}(z^{\sigma+n+1})
\quad \textup{quand}\quad z\to 0 \,;\\
\dis I(z)=\sum_{j=1}^A
Q_{j}(z)\frac{\log^{j-1}(1/z)}{(j-1)!}=\mathcal{O}(\text{?})
\quad \textup{quand}\quad z\to 1.
\end{cases}
\end{equation}
Il nous reste \`a montrer que l'on peut
mettre $\mathcal{O}((z-1)^{A(n+1)-\rho-\sigma-2})$ \`a la place  de
$\mathcal{O}(\text{?})$. Pour le faire, supposons donn{\'e}e
une fonction  $f(z;q)$ analytique en $z$ dans un voisinage ouvert
suffisamment grand de 1, et telle que
\begin{equation} \label{eq:f}
f(z;q)=\mathcal{O}(z-q^{-\ell})\quad
\text{quand}\quad z\to q^{-\ell}
\end{equation}
pour tout
$\ell\in\{-m,-m+1,\dots,p\}$. Pour simplifier, on consid{\`e}re
le cas o{\`u} $p=0$, car l'argument qui va suivre  
se g{\'e}n{\'e}ralise sans difficult\'e au
cas o{\`u} $p$ est quelconque. 
On d{\'e}veloppe tout d'abord  $f(z;q)$ en s{\'e}rie
de Taylor autour de $z=1$~:
\begin{equation} \label{eq:Taylor}
f(z;q)=
\sum _{k=1} ^{\infty}f_k(q)(z-1)^k,
\end{equation}
ce d{\'e}veloppement ne contenant pas de terme constant {\`a}
cause de la condition \eqref{eq:f} pour $\ell=0$.
Nous supposons alors une condition suppl{\'e}mentaire~:
la limite $
\lim _{q\to1} ^{}f_k(q)$ existe pour tout $k$ et
$$
\lim _{q\to1} ^{}\sum _{k=1} ^{\infty}f_k(q)(z-1)^k=
\sum _{k=1} ^{\infty}\lim _{q\to1} ^{}f_k(q)(z-1)^k.$$
Cette condition est satisfaite
dans notre cas, c'est-{\`a}-dire pour
\begin{equation} \label{eq:special}
f(z;q)=
-(1-q)^{A(n+1)-\sigma-\rho-1}\sum_{j=1}^A
P_{j}(zq^{1-j};q)\frac{\(-\log_q(1/z)\)_{j-1}}{(j-1)!}.
\end{equation}

La condition
\eqref{eq:f} pour les valeurs non-nulles de $\ell$ implique le
syst{\`e}me d'{\'e}quations
\begin{equation} \label{eq:fk}
0=f(q^{-\ell};q)=
\sum _{k=1} ^{\infty}f_k(q)(q^{-\ell}-1)^k,\quad \ell\in\{-m,-m+1,\dots,-1\}.
\end{equation}
On multiplie la $\ell$-i{\`e}me {\'e}quation par
$$c_\ell=\frac {1} {(q^{-\ell}-1)
\underset{h\ne -\ell}{\prod\limits _{h=1} ^{m}}(q^{-\ell}-q^h)}$$
et en faisant la somme de ces {\'e}quations multipli{\'ees} par
le facteur $c_\ell$ correspondant
sur $\ell\in\{-m,-m+1,\dots, -1\}$, on obtient
\begin{align} \notag
0&=\sum _{\ell=-m}
^{-1}c_\ell f(q^{-\ell};q)\\
\notag
&=\sum _{\ell=-m}
^{-1}c_\ell \sum _{k=1} ^{\infty}f_k(q)(q^{-\ell}-1)^k\\
\label{eq:fkq}
&=\sum _{k=1} ^{\infty}f_k(q)\sum _{\ell=-m}
^{-1}c_\ell (q^{-\ell}-1)^k .
\end{align}
{\`A} ce point, on note que le choix des coefficients
$c_\ell$ implique
que les coefficients de $f_k(q)$ dans la somme \eqref{eq:fkq}
sont nuls pour $k=1,2,\dots,m-1$, et que le coefficient de
$f_m(q)$ est exactement~1. Dans le r{\'e}sultat
$$0=f_m(q)+\sum _{k\ge m+1} ^{\infty}f_k(q)\sum _{\ell=-m}
^{-1}\frac {(q^{-\ell}-1)^k} {(q^{-\ell}-1)
\underset{h\ne -\ell}{\prod\limits _{h=1} ^{m}}(q^{-\ell}-q^h)}, 
$$
on fait maintenant tendre $q$ vers 1~: comme la somme ext{\'e}rieure porte sur
les $k>m$, la limite du sommande
de la somme int{\'e}rieure est toujours z{\'e}ro.
Par cons{\'e}quent, on obtient bien que $f_m(1)=0$. De plus,
si l'on applique le m{\^e}me argument pour $\ell\in\{-\bar m,
-\bar m+1,\dots,-1\}$ avec $\bar m=m-1, m-2,\dots,1$, alors
on obtient que $f_{\bar m}(1)=0$ pour tout $\bar m\in
\{1,2,\dots,m\}$, ce qui prouve que
$$
f(z;1)=\mathcal O\big((z-1)^m\big).
$$

Cet argument, appliqu{\'e} {\`a} \eqref{eq:special}, montre
que l'on peut bien remplacer $\mathcal O(\text{?})$ dans
\eqref{problemepadeQ} par
$\mathcal{O}((z-1)^{A(n+1)-\rho-\sigma-2})$, comme annonc\'e.
Le probl{\`e}me d'approximation \eqref{problemepadeQ} est
donc exactement le probl{\`e}me \eqref{problemepade}.
Comme il est d{\'e}montr{\'e} dans \cite{firi} que ce probl{\`e}me a
une solution unique, on a forc{\'e}ment l'{\'e}galit{\'e}
$Q_j(z)=P_j(z)$ pour tout $j$.

\def\refname{Bibliographie}

\end{document}